\documentclass[11pt,reqno]{amsart}
\usepackage{amsmath}
\usepackage{amssymb}
\usepackage{enumerate}
\usepackage[bookmarks]{hyperref}

\newtheorem{theorem}{Theorem}[section]
\newtheorem{lemma}[theorem]{Lemma}
\newtheorem{corollary}[theorem]{Corollary}
\newtheorem{proposition}[theorem]{Proposition}

\theoremstyle{definition}

\numberwithin{equation}{section}

\newcommand{\ind}{{\rm ind}}

\allowdisplaybreaks[1]

\begin{document}

\title[Toeplitz operators on harmonic Bergman spaces]{Toeplitz operators\\ on generalized harmonic Bergman spaces}
\author{Trieu Le}
\address{Department of Mathematics, Mail Stop 942, University of Toledo, Toledo, OH 43606}
\email{trieu.le2@utoledo.edu}


\begin{abstract}
We study Toeplitz operators with uniformly continuous symbols on generalized harmonic Bergman spaces of the unit ball in $\mathbb{R}^n$. We describe their essential spectra and establish a short exact sequence associated with the $C^{*}$-algebra generated by these operators.
\end{abstract}
\maketitle

\section{Introduction}\label{S:Intro}

Let $n\geq 2$ be a fixed integer. We write $\mathbb{B}$ for the open unit ball and $\mathbb{S}$ for the unit sphere in $\mathbb{R}^n$. The closure of $\mathbb{B}$, which is the closed unit ball, is denoted by $\bar{\mathbb{B}}$. For any $x=(x_1,\ldots, x_n)$ in $\mathbb{R}^n$, we use $|x|$ to denote the Euclidean norm of $x$, that is, $|x|=(x_1^2+\cdots+x_n^2)^{1/2}$.

Let $\nu$ be a regular Borel probability measure on $\mathbb{B}$ that is invariant under the action of the group of orthogonal transformations $O(n)$. Then there is a regular Borel probability measure $\mu$ on the interval $[0,1)$ so that the integration in polar coordinates formula $$\int_{\mathbb{B}}f(x)\mathrm{d}\nu(x)=\int_{[0,1)}\int_{\mathbb{S}}f(r\zeta)\mathrm{d}\sigma(\zeta)\mathrm{d}\mu(r)$$
holds for all functions $f$ that belong to $L^1(\mathbb{B},\nu)$. Here $\sigma$ is the unique $O(n)$-invariant regular Borel probability measure on the unit sphere $\mathbb{S}$. We are interested in measures $\nu$ whose support is not entirely contained in a compact subset of the unit ball so we will assume through out the paper that $\nu(\{x\in\mathbb{B}: |x|\geq r\})>0$ for all $0<r<1$. This is equivalent to the condition that $\mu([r,1))>0$ for all $0<r<1$. 

The (generalized) harmonic Bergman space $b^2_{\nu}$ is the space of all harmonic functions that belong also to the Hilbert space $L^2_{\nu}=L^2(\mathbb{B},\nu)$. It follows from Poisson integral representation of harmonic functions and the assumption about $\nu$ that for any compact subset $K$ of $\mathbb{B}$, there is a constant $C_{K}$ such that
\begin{equation}\label{Eqn:kernel}
|u(x)|\leq C_K\|u\|=\Big(\int_{\mathbb{B}}|u(x)|^2\mathrm{d}\nu(x)\Big)^{1/2}
\end{equation}
for all $x$ in $K$ and all $u$ in $b^2_{\nu}$. This implies that $b^2_{\nu}$ is a closed subspace of $L^2_{\nu}$ and that the evaluation map $u\mapsto u(x)$ is a bounded linear functional on $b^2_{\nu}$ for each $x$ in $\mathbb{B}$. By the Riesz's representation, there is a function $R_x$ in $b^2_{\nu}$ so that $u(x)=\langle u, R_x\rangle$. The function $R(y,x):=R_x(y)$ for $x,y\in\mathbb{B}$ is called the reproducing kernel for $b^2_{\nu}$.

Let $Q$ denote the orthogonal projection from $L^2_{\nu}$ onto $b^2_{\nu}$. For a bounded measurable function $f$ on $\mathbb{B}$, the Toeplitz operator $T_{f}: b^2_{\nu}\to b^2_{\nu}$ is defined by $$T_fu = QM_fu = Q(fu), \quad u\in b^2_{\nu}.$$
Here $M_f: L^2_{\nu}\to L^2_{\nu}$ is the operator of multiplication by $f$. The function $f$ is called the symbol of $T_f$. We also define the Hankel operator $H_f: b^2_{\nu}\to (b^2_{\nu})^{\perp}$ by $$H_fu = (1-Q)M_fu = (1-Q)(fu),\quad u\in b^2_{\nu}.$$ It is immediate that $\|T_f\|\leq\|f\|_{\infty}$ and $\|H_{f}\|\leq\|f\|_{\infty}$.

For $f,g$ bounded measurable functions on $\mathbb{B}$, the following basic properties are immediate from the definition of Toeplitz and Hankel operators: 
\begin{equation}\label{Eqn:ToeplitzHankel}
T_{gf}-T_gT_f = H^{*}_{\bar{g}}H_{f},
\end{equation} and $$(T_g)^{*}=T_{\bar{g}},\quad\quad T_{af+bg}=aT_f+bT_g,$$
where $a, b$ are complex numbers and $\bar{g}$ denotes the complex conjugate of $g$.

If $\mathrm{d}\nu(x)=\mathrm{d}V(x)$, where $V$ is the normalized Lebesgue volume measure on $\mathbb{B}$, then $b^2_{\nu}$ is the usual unweighted harmonic Bergman space. See \cite[Chapter 8]{AxlerSpringer2001} for more details about this space. If $\mathrm{d}\nu(x)=c_{\alpha}(1-|x|^2)^{\alpha}\mathrm{d}V(x)$ where $-1<\alpha<\infty$ and $c_{\alpha}$ is a normalizing constant, then $\nu$ is a weighted Lebesgue measure on $\mathbb{B}$ and $b^2_{\nu}$ is a weighted harmonic Bergman space.  Compactness of certain classes of  Toeplitz operators on these weighted harmonic Bergman spaces was considered by K. Stroethoff in \cite{StroethoffAuMS1998}. He also described the essential spectra of Toeplitz operators with uniformly continuous symbols. He showed that if $f$ is a continuous function on the closed unit ball $\bar{\mathbb{B}}$, then the essential spectrum of $T_f$ is the same as the set $f(\mathbb{S})$. This result in the setting of unweighted harmonic Bergman spaces was obtained earlier by J. Miao \cite{MiaoIEOT1997}. More recently, B.R. Choe, Y.J. Lee and K. Na \cite{ChoeNMJ2004} showed that the above essential spectral formula remains valid for unweighted harmonic Bergman space of any bounded domain with smooth boundary in $\mathbb{R}^n$. The common approach, which was used in all of the above papers, involved a careful estimate on the kernel function. In the case where $\nu$ is not a weighted Lebesgue measure, it seems that similar estimates are not available. Nevertheless, with a different approach, we still obtain the aforementioned essential spectral formula.

Let $\mathfrak{T}$ denote the $C^{*}$-algebra generated by all Toeplitz operators $T_f$, where $f$ belongs to the space $C(\bar{\mathbb{B}})$ of continuous functions on the closed unit ball. Let $\mathfrak{CT}$ denote the two-sided ideal of $\mathfrak{T}$ generated by commutators $[T_f,T_g]=T_fT_g-T_gT_f$, for $f,g\in C(\bar{\mathbb{B}})$. In the case $n=2$ and $\nu$ the normalized Lebesgue measure on the unit disk, K. Guo and D. Zheng \cite{GuoJMAA2002} proved that $\mathfrak{CT}=\mathcal{K}$, the ideal of compact operators on $b^2_{\nu}$, and there is a short exact sequence $$0\rightarrow\mathcal{K}\rightarrow\mathfrak{T}\rightarrow C(\mathbb{S})\rightarrow 0.$$ They also proved that any Fredholm operator in the Toeplitz algebra $\mathfrak{T}$ has Fredholm index $0$. We will show that these results are in fact valid for all $n\geq 2$.

The paper is organized as follows. In Section 2 we give some preliminaries. We then study Toeplitz operators with uniformly continuous symbols and establish the essential spectral formula in Section 3. The Toeplitz algebra and the associated short exact sequence are studied in Section 4. We close the paper with a criterion for compactness of operators with more general symbols in Section 5.

\section{Preliminaries}

It is well known that the reproducing kernel $R(x,y)$ is symmetric and real-valued for $x,y\in\mathbb{B}$. From \eqref{Eqn:kernel} we see that for any compact subset $K$ of $\mathbb{B}$ and $x\in K$,
\begin{equation}\label{Eqn:ineqkernel}
R(x,x) = R_x(x)\leq C_K\|R_x\|=C_K(\langle R_x,R_x\rangle)^{1/2} = C_K(R(x,x))^{1/2}.
\end{equation}
This shows that $0\leq R(x,x)\leq C_{K}^2$ for $x\in K$. So the function $x\mapsto K(x,x)$ is bounded on compact subsets of $\mathbb{B}$.

A polynomial $p$ in the variable $x$ with complex coefficients is homogeneous of degree $m$ (or $m$-homogeneous), where $m\geq 0$ is an integer, if $p(tx)=t^mp(x)$ for all non-zero real numbers $t$. We write $\mathcal{P}_{m}$ for the vector space of all $m$-homogeneous polynomials on $\mathbb{R}^n$. We use $\mathcal{H}_m$ to denote the subspace of $\mathcal{P}_m$ consisting of harmonic polynomials. The subspace $\mathcal{H}_m$ is finite dimensional and its dimension $h_m$ is given by $h_0=1, h_1=n$ and $h_m=\binom{n+m-1}{n-1}-\binom{n+m-3}{n-1}$ for $m\geq 2$. See \cite[Proposition 5.8]{AxlerSpringer2001}.

For polynomials $p$ in $\mathcal{H}_m$ and $q$ in $\mathcal{H}_k$, using the orthogonality of their restrictions on the sphere \cite[Proposition 5.9]{AxlerSpringer2001} and integration in polar coordinates we obtain
\begin{align}\label{Eqn:orthogonality}
\langle p, q\rangle & = \Big(\int_{[0,1)}r^{m+k}\mathrm{d}\mu(r)\Big)\int_{\mathbb{S}}p\bar{q}\mathrm{d}\sigma\notag\\
& = \begin{cases}
		0 & \text{ if } m\neq k\\
		\Big(\int_{[0,1)}r^{2m}\mathrm{d}\mu(r)\Big)\int_{\mathbb{S}}p\bar{q}\mathrm{d}\sigma & \text{ if } m=k.
    \end{cases}
\end{align}
This shows that the spaces $\mathcal{H}_m$ for $m=0,1,\ldots$ are pairwise orthogonal. On the other hand, \cite[Corollary 5.34]{AxlerSpringer2001} shows that if $u$ is a harmonic function on the unit ball $\mathbb{B}$, then there exist polynomials $p_m\in\mathcal{H}_m$ such that $u(x)=\sum_{m=0}^{\infty}p_m(x)$ for all $x$ in $\mathbb{B}$. The series converges uniformly on compact subsets of $\mathbb{B}$. Thus we have the orthogonal decomposition $b^2_{\nu} = \oplus_{m=0}^{\infty}\mathcal{H}_m$.

We next present some other elementary results which we will use later in the paper. The following lemmas follow from the above orthogonal decomposition of $b^2_{\nu}$. Since these are well known results in functional analysis, we omit their proofs.

\begin{lemma}
\label{L:cpt}
Suppose $A$ is a compact operator from $b^2_{\nu}$ into a Hilbert space $\mathcal{L}$. Then $\lim_{m\to\infty}\|A|_{\mathcal{H}_m}\|=0$.
\end{lemma}

The converse of Lemma \ref{L:cpt} is false in general. However, if additional conditions are imposed on the images of the subspaces $\mathcal{H}_m$ under $A$, the converse holds.

\begin{lemma}
\label{L:cptOperator}
Suppose $A$ is an operator defined on the algebraic direct sum of the subspaces $\mathcal{H}_m$ into a Hilbert space $\mathcal{L}$ so that $A(\mathcal{H}_m)\perp A(\mathcal{H}_l)$ for all $m\neq l$ and $\lim_{m\to\infty}\|A|_{\mathcal{H}_m}\| = 0$. Then $A$ extends (uniquely) to a compact operator from $b^2_{\nu}$ into $\mathcal{L}$.
\end{lemma}

We will also need the following lemma, which is a special case of \cite[Lemma 2.4]{LePAMS2009}.
\begin{lemma}\label{L:limitIntegrals}
Suppose $\varphi$ is a function on $[0,1)$ so that $\lim_{r\uparrow 1}\varphi(r)=\gamma$, then $$\lim_{m\to\infty}\dfrac{\int_{[0,1)}\varphi(r)r^{2m}\mathrm{d}\mu(r)}{\int_{[0,1)}r^{2m}\mathrm{d}\mu(r)}=\gamma.$$
\end{lemma}

\section{Toeplitz operators with uniformly continuous symbols}

In this section we study Toeplitz operators whose symbols are continuous on the closed unit ball. We show that the essential spectrum of such an operator is the set of values of its symbol on the unit sphere.

For any integer $k\geq 0$, we denote the $2k$-th moment of $\mu$ by $\hat{\mu}(k)$, that is, $$\hat{\mu}(k)=\int_{[0,1)}r^{2k}\mathrm{d}\mu(r) = \int_{\mathbb{B}}|x|^{2k}\mathrm{d}\nu(x).$$
For any $u\in b^2_{\nu}$, write $u=\sum_{m=0}^{\infty}u_m$, where $u_m\in\mathcal{H}_m$ for $m\geq 0$. We then have
\begin{align*}
\|u\|^2 & = \sum_{m=0}^{\infty}\|u_m\|^2 = \sum_{m=0}^{\infty}\int_{\mathbb{B}}|u_m(x)|^2\mathrm{d}\nu(x)\\
& = \sum_{m=0}^{\infty}\int_{[0,1)}r^{2m}\mathrm{d}\mu(r)\int_{\mathbb{S}}|u_m(\zeta)|^2\mathrm{d}\sigma(\zeta)\\
& = \sum_{m=0}^{\infty}\hat{\mu}(m)\int_{\mathbb{S}}|u_m(\zeta)|^2\mathrm{d}\sigma(\zeta).
\end{align*}
This shows that the linear map $W: b^2_{\nu}\longrightarrow L^2(\mathbb{S})$ defined by $$W(u) = \sum_{m=0}^{\infty}(\hat{\mu}(m))^{1/2}\ u_m|_{\mathbb{S}}$$ is isometric.

The restriction of an element in $\mathcal{H}_m$ to $\mathbb{S}$ is called a spherical harmonic of degree $m$. Theorem 5.12 in \cite{AxlerSpringer2001} shows that the span of all spherical harmonics is dense in $L^2(\mathbb{S})$. We then conclude that $W$ is a surjective isometry, hence a unitary operator.

For a continuous function $f$ on the closed unit ball, let $f_{*}$ denote the restriction of $f$ on the unit sphere. Recall that the operator $M_f$ is the multiplication operator on $L^2_{\nu}$ with symbol $f$. As usual, we denote its restriction on $b^2_{\nu}$ by $M_f|_{b^2_{\nu}}$. We also write $M_{f_{*}}$ for the multiplication operator on $L^2(\mathbb{S})$ with symbol $f_{*}$. Using the above unitary, we establish the following connection between these two operators.

\begin{theorem}\label{T:reduction}
Let $f$ be in $C(\bar{\mathbb{B}})$. Then the operator $$M_f|_{b^2_{\nu}}-W^{*}M_{f_{*}}W: b^2_{\nu}\longrightarrow L^2_{\nu}$$ is compact.
\end{theorem}

\begin{proof}
Let $$\mathcal{A}=\{f\in C(\bar{\mathbb{B}}): M_{f}|_{b^2_{\nu}}-W^{*}M_{f_{*}}W \text{ is compact}\}.$$ We need to show that $\mathcal{A}=C(\bar{\mathbb{B}})$. It is clear that $\mathcal{A}$ is a closed linear subspace of $C(\bar{\mathbb{B}})$. Now suppose $f,g$ are in $\mathcal{A}$. Then there are compact operators $K_f$ and $K_g$ from $b^2_{\nu}$ into $L^2_{\nu}$ so that $$M_{f}|_{b^2_{\nu}}=W^{*}M_{f_{*}}W+K_f\text{ and } M_{g}|_{b^2_{\nu}}=W^{*}M_{g_{*}}W+K_g.$$
Since the range of $W^{*}$ is $b^2_{\nu}$, we have $(I-Q)W^{*}M_{g_{*}}W=0$ (we recall here that $Q$ is orthogonal projection from $L^2_{\nu}$ onto $b^2_{\nu}$). This implies that $(1-Q)M_{g}|_{b^2_{\nu}}=(1-Q)K_g$ and so we have
\begin{align*}
M_{fg}|_{b^2_{\nu}} & = M_fM_g|_{b^2_{\nu}} = M_f(1-Q)M_g|_{b^2_{\nu}}+M_fQM_g|_{b^2_{\nu}}\\
& = M_f(1-Q)K_g+(W^{*}M_{f_{*}}W+K_f)(W^{*}M_{g_{*}}W+QK_g)\\
& = W^{*}M_{f_{*}}WW^{*}M_{g_{*}}W+K\\
& = W^{*}M_{(fg)_{*}}W+K,
\end{align*}
where $K$ is compact. Thus, $fg\in\mathcal{A}$ if $f,g$ are in $\mathcal{A}$. We have showed that $\mathcal{A}$ is a closed subalgebra of $C(\bar{\mathbb{B}})$. To complete the proof of the theorem, we will show that the coordinate functions $x_1,\ldots,x_n$ belong to $\mathcal{A}$. By symmetry, we only need to check that $f(x)=x_1$ is in $\mathcal{A}$.

For any integer $m\geq 0$ and $p\in\mathcal{H}_m$, the polynomial $fp$ is homogeneous of degree $m+1$. Since $\Delta^{2}(fp)=\Delta^{2}(x_1p(x))=0$, \cite[Theorem 5.21]{AxlerSpringer2001} shows that there is a unique decomposition $$M_{f}p = fp = p_{m+1}+|x|^2p_{m-1},$$
where $p_{m+1}\in\mathcal{H}_{m+1}$ and $p_{m-1}\in\mathcal{H}_{m-1}$ if $m\geq 1$ and $p_{m-1}=0$ if $m=0$. Using integration in polar coordinates and the fact that restrictions of homogeneous harmonic polynomials of different degrees are orthogonal in $L^2(\mathbb{S})$, we obtain 
\begin{equation}\label{Eqn:multiplication}
\|M_{f}p\|^2 = \|p_{m+1}\|^2 + \||x|^2p_{m-1}\|^2 = \|p_{m+1}\|^2 + \dfrac{\hat{\mu}(m+1)}{\hat{\mu}(m-1)}\|p_{m-1}\|^2.
\end{equation}
On the other hand,
\begin{align*}
W^{*}M_{f_{*}}W(p) & = (\hat{\mu}(m))^{1/2}W^{*}M_{f_{*}}(p|_{\mathbb{S}})\\
& = (\hat{\mu}(m))^{1/2}W^{*}(fp|_{\mathbb{S}})\\
& = (\hat{\mu}(m))^{1/2}W^{*}(p_{m+1}|_{\mathbb{S}}+p_{m-1}|_{\mathbb{S}})\Big)\\
& = \Big(\dfrac{\hat{\mu}(m)}{\hat{\mu}(m+1)}\Big)^{1/2}p_{m+1}+\Big(\dfrac{\hat{\mu}(m)}{\hat{\mu}(m-1)}\Big)^{1/2}p_{m-1}
\end{align*}

Therefore,
\begin{align}\label{Eqn:sumCptOperators}
& (M_f|_{b^2_{\nu}}-W^{*}M_{f_{*}}W)(p)\\
&\quad\quad\quad = \Big\{1-\Big(\dfrac{\hat{\mu}(m)}{\hat{\mu}(m+1)}\Big)^{1/2}\Big\}p_{m+1}+\Big(|x|^2-\Big(\dfrac{\hat{\mu}(m)}{\hat{\mu}(m-1)}\Big)^{1/2}\Big)p_{m-1}.\notag
\end{align}

Now we define
\begin{align*}
A_1(p) & = \Big\{1-\Big(\dfrac{\hat{\mu}(m)}{\hat{\mu}(m+1)}\Big)^{1/2}\Big\}p_{m+1}, \text{ and }\\
A_2(p) & = \Big\{|x|^2-\Big(\dfrac{\hat{\mu}(m)}{\hat{\mu}(m-1)}\Big)^{1/2}\Big\}p_{m-1},
\end{align*}
for $p\in\mathcal{H}_m$ and extend $A_1$ and $A_2$ by linearity to the algebraic direct sum of the subspaces $\mathcal{H}_m$, $m=0,1\ldots$. Using \eqref{Eqn:multiplication} and integration in polar coordinates, we have
\begin{align*}
\|A_1(p)\| & = \Big|1-\Big(\dfrac{\hat{\mu}(m)}{\hat{\mu}(m+1)}\Big)^{1/2}\Big|\|p_{m+1}\|\leq \Big\{\Big(\dfrac{\hat{\mu}(m)}{\hat{\mu}(m+1)}\Big)^{1/2}-1\Big\}\|p\|,
\end{align*}
\begin{align*}
&\|A_2(p)\|^2\\
& = \int_{[0,1)}\Big(r^2-\big(\dfrac{\hat{\mu}(m)}{\hat{\mu}(m-1)}\big)^{1/2}\Big)^2 r^{2m-2}\mathrm{d}r\int_{\mathbb{S}}|p_{m-1}(\zeta)|^2\mathrm{d}\sigma(\zeta)\\
& = \Big\{\hat{\mu}(m+1)+\hat{\mu}(m)-2\hat{\mu}(m)\big(\dfrac{\hat{\mu}(m)}{\hat{\mu}(m-1)}\big)^{1/2}\Big\}\big(\hat{\mu}(m-1)\big)^{-1}\|p_{m-1}\|^2\\
& \leq \Big\{\hat{\mu}(m+1)+\hat{\mu}(m)-2\hat{\mu}(m)\big(\dfrac{\hat{\mu}(m)}{\hat{\mu}(m-1)}\big)^{1/2}\Big\}\big(\hat{\mu}(m+1)\big)^{-1}\|p\|^2.
\end{align*}
From Lemma \ref{L:limitIntegrals} we have $\lim_{k\to\infty}\hat{\mu}(k)/\hat{\mu}(k+1) = 0$. This implies that $\lim_{m\to\infty}\|A_j|_{\mathcal{H}_m}\|=0$ for $j=1,2$. On the other hand, by the orthogonality of homogeneous harmonic polynomials of different degrees when restricted to the sphere, we see that $A_j(\mathcal{H}_m)\bot A_j(\mathcal{H}_k)$ if $m\neq k$ for $j=1,2$. Lemma \ref{L:cptOperator} now shows that $A_1$ and $A_2$ extend to compact operators on $b^2_{\nu}$. From \eqref{Eqn:sumCptOperators}, $M_f|_{b^2_{\nu}}-W^{*}M_{f_{*}}W$, being a sum of two compact operators, is also compact. This completes the proof of the theorem.
\end{proof}

Theorem \ref{T:reduction} has important consequences that we now describe. Since the image of $W^{*}$ is contained in $b^2_{\nu}$, we have $QW^{*}=W^{*}$ and $(1-Q)W^{*}=0$. Theorem \ref{T:reduction} implies that the Toeplitz operator $T_{f}=QM_{f}|_{b^2_{\nu}}$ is a compact perturbation of $W^{*}M_{f_{*}}W$ and the Hankel operator $H_{f}=(1-Q)M_{f}|_{b^2_{\nu}}$ is compact for any $f$ in $C(\bar{\mathbb{B}})$. Now for any bounded measurable function $g$ on $\mathbb{B}$, \eqref{Eqn:ToeplitzHankel} shows that both operators $T_{gf}-T_gT_f$ and $T_{gf}-T_{f}T_g$ are compact.

Let us write $\mathcal{B}(b^2_{\nu})$ for the $C^{*}$-algebra of all bounded operators on $b^2_{\nu}$. Let $\mathcal{K}$ denote the ideal of all compact operators on $b^2_{\nu}$. For any bounded operator $A$, recall that the essential spectrum of $A$, denoted by $\sigma_{e}(A)$, is the spectrum of $A+\mathcal{K}$ in the quotient algebra $\mathcal{B}(b^2_{\nu})/\mathcal{K}$. The essential norm $\|A\|_{e}$ is the norm of $A+\mathcal{K}$ in $\mathcal{B}(b^2_{\nu})/\mathcal{K}$.

If $f_{*}$ is the restriction of $f$ on the unit sphere $\mathbb{S}$, then $\sigma_{e}(M_{f_{*}})=f_{*}(\mathbb{S})=f(\mathbb{S})$ and $\|M_{f_{*}}\|_{e}=\sup\{|f_{*}(\zeta)|: \zeta\in\mathbb{S}\}=\sup\{|f(\zeta)|: \zeta\in\mathbb{S}\}$. Theorem \ref{T:reduction} now implies the following results, which were obtained earlier by Stroethoff \cite{StroethoffGMJ1997,StroethoffAuMS1998} for weighted spaces with a different approach.

\begin{corollary}\label{C:essSpectra} Let $f$ be a uniformly continuous function and $g$ be a bounded measurable function on $\mathbb{B}$. Then
$T_{gf}-T_gT_f$ and $T_{gf}-T_fT_g$ are compact operators, $\sigma_{e}(T_f)=f(\mathbb{S})$ and $\|T_f\|_{e}=\sup\{|f(\zeta)|: \zeta\in\mathbb{S}\}$. In particular, $T_f$ is compact if and only if $f$ vanishes on $\mathbb{S}$. 
\end{corollary}

\section{The Toeplitz Algebra}

We now turn our attention to the $C^{*}$-algebra $\mathfrak{T}$ generated by all Toeplitz operators $T_f$, where $f$ belongs to $C(\bar{\mathbb{B}})$. Our main result in this section is a description of this algebra as an extension of the compact operators $\mathcal{K}$ by continuous functions on the unit sphere. We begin by exhibiting a class of block diagonal operators in $\mathfrak{T}$.

A function $f$ on the unit ball is called radial if there is a function $\varphi$ defined on the interval $[0,1)$ so that $f(x)=\varphi(|x|)$ for $\nu$-almost every $x$ in $\mathbb{B}$. The following lemma, which is Lemma 4.2 in \cite{StroethoffAuMS1998} in the case $\nu$ a weighted Lebesgue measure, shows that each $\mathcal{H}_m, m=0,1,\ldots$ is an eigenspace for $T_f$. For completeness we include here a proof.

\begin{lemma}
\label{L:radial}
If $f$ is a bounded radial function on $\mathbb{B}$, then each non-zero homogeneous harmonic polynomial of degree $m\geq 0$ is an eigenvector of $T_f$ with eigenvalue given by
\begin{equation*}
	\lambda_m = \dfrac{\int_{[0,1)}\varphi(r)r^{2m}\mathrm{d}\mu(r)}{\int_{[0,1)}r^{2m}\mathrm{d}\mu(r)},
\end{equation*}
where $\varphi$ is a bounded function on the interval $[0,1)$ so that $f(x)=\varphi(|x|)$ for $\nu$-almost every $x\in\mathbb{B}$.
\end{lemma}

\begin{proof}
For any homogeneous harmonic polynomials $p$ of degree $m$ and $q$ of degree $k$, using \eqref{Eqn:orthogonality} we have
		\begin{align*}
			\langle T_fp,q\rangle & = \langle fp,q\rangle = \int_{[0,1)}\int_{\mathbb{S}}f(r\zeta)p(r\zeta)\bar{q}(r\zeta)\mathrm{d}\sigma(\zeta)\mathrm{d}\mu(r)\\
								  & = \Big(\int_{[0,1)}\varphi(r)r^{m+k}\mathrm{d}\mu(r)\Big)\int_{\mathbb{S}}p\bar{q}\mathrm{d}\sigma\\
								  & = \begin{cases}
								    0 & \text{ if } m\neq k,\\
									\Big(\int_{[0,1)}\varphi(r)r^{2m}\mathrm{d}\mu(r)\Big)\int_{\mathbb{S}}p\bar{q}\mathrm{d}\sigma & \text{ if } m=k
								  \end{cases}\\
								  & =\lambda_{m}\langle p, q\rangle.
		\end{align*}
Since the span of homogeneous harmonic polynomials is dense in $b^2_{\nu}$, we conclude that $T_fp = \lambda_{m}p$ for any $p$ in $\mathcal{H}_m$.
\end{proof}

Let $\eta(x)=|x|^2$ for $x$ in $\mathbb{B}$. From the lemma, each subspace $\mathcal{H}_m$ is an eigenspace for $T_{\eta}$ with corresponding eigenvalue $\gamma_{m}=\dfrac{\int_{[0,1)}{r^{2m+2}\mathrm{d}\mu(r)}}{\int_{[0,1)}r^{2m}\mathrm{d}\mu(r)}$. We show that these eigenvalues form a strictly increasing sequence. In particular, they are pairwise distinct.

\begin{lemma}
\label{L:inequality}
For any integer $m\geq 0$ we have $$\dfrac{\int_{[0,1)}{r^{2m+2}\mathrm{d}\mu(r)}}{\int_{[0,1)}r^{2m}\mathrm{d}\mu(r)} < \dfrac{\int_{[0,1)}{r^{2(m+1)+2}\mathrm{d}\mu(r)}}{\int_{[0,1)}r^{2(m+1)}\mathrm{d}\mu(r)}.$$
\end{lemma}

\begin{proof}
	Let $a(r)=r^{m}$ and $b(r)=r^{m+2}$ for $0\leq r<1$. Cauchy-Schwarz's inequality gives $$\Big(\int a(r)b(r)\mathrm{d}\mu(r)\Big)^2\leq\Big(\int_{[0,1)}a^2(r)\mathrm{d}\mu(r)\Big)\Big(\int_{[0,1)}b^2(r)\mathrm{d}\mu(r)\Big),$$
	which is almost the same as the required inequality. We need to show why the equality cannot occur. This is indeed the case because the ratio $b(r)/a(r)=r^2$ is not a constant function $\mu$-almost everywhere on $(0,1)$.
\end{proof}

We are now ready for the description of $\mathfrak{T}$, which is the main result in this section.

\begin{theorem}\label{T:ToeplitzAlgebra}
The following statements hold.

(1) The commutator ideal $\mathfrak{CT}$ of $\mathfrak{T}$ is the same as the ideal $\mathcal{K}$ of compact operators on $b^2_{\nu}$. 

(2) Any element of $\mathfrak{T}$ has the form $T_f+K$ for some $f$ in $C(\bar{\mathbb{B}})$ and $K\in\mathcal{K}$ and there is a short exact sequence $$0\rightarrow\mathcal{K}\rightarrow\mathfrak{T}\rightarrow C(\mathbb{S})\rightarrow 0.$$
\end{theorem}

Guo and Zheng \cite{GuoJMAA2002} proved Theorem \ref{T:ToeplitzAlgebra} for the case $n=2$ and $\nu$ the normalized Lebesgue measure. 

When $n$ is an even number and $\nu$ is the normalized Lebesgue measure, Theorem \ref{T:ToeplitzAlgebra} was proved by L. Coburn \cite{CoburnIMJ1973} for Toeplitz operators on the holomorphic Bergman space. The idea of our proof is similar to that of \cite[Theorem 1]{CoburnIMJ1973}.

\begin{proof}
We first show that $\mathfrak{T}$ is irreducible on $b^2_{\nu}$. Suppose $A$ is an operator on $b^2_{\nu}$ that commutes with all elements of $\mathfrak{T}$. Then in particular, $A$ commutes with $T_{\eta}$. For any homogeneous polynomials $p$ of degree $m$ and $q$ of degree $k$, we have
\begin{align*}
	\langle AT_{\eta}p,q\rangle & = \gamma_m\langle Ap,q\rangle,\\
	\langle T_{\eta}Ap,q\rangle & = \langle Ap, T_{\bar{\eta}}q\rangle = \langle Ap, T_{\eta}q\rangle=\gamma_k\langle Ap,q\rangle.
\end{align*}
Since $AT_{\eta}=T_{\eta}A$ and $\gamma_m\neq\gamma_k$ if $m\neq k$, we conclude that $\langle Ap,q\rangle=0$ if $m\neq k$. This implies that each subspace $\mathcal{H}_m$ is invariant, hence also reducing for $A$. In particular, $\mathcal{H}_{0}$ reduces $A$. But $\mathcal{H}_{0}$ is a one-dimensional space spanned by the constant function $e_{0}(x)=1$, so we have $Ae_{0}=\lambda e_{0}$ for some scalar $\lambda$.

For each harmonic polynomial $p$, we have $$A(p) = AT_{p}(e_{0}) = T_{p}A(e_{0})=\lambda T_{p}(e_{0}) = \lambda p.$$
Since the space of harmonic polynomials is dense in $b^2_{\nu}$, we conclude that $A = \lambda I_{b^2_{\nu}}$. Thus $\mathfrak{T}$ is irreducible. Now Corollary \ref{C:essSpectra} shows that $\mathfrak{T}$ contains a non-zero compact operator (for example $T_{1-|x|^2}$). It then follows from a well known result in $C^{*}$-algebra \cite[Theorem 5.39]{Douglas1972} that $\mathfrak{T}$ contains the ideal $\mathcal{K}$ of compact operators. Therefore the commutator ideal $\mathfrak{CT}$ contains the commutator ideal of $\mathcal{K}$, which is the same as $\mathcal{K}$.

On the other hand, for any functions $f, g$ in $C(\bar{\mathbb{B}})$, the commutator $T_fT_g - T_g T_f = (T_fT_g-T_{fg})-(T_gT_f-T_{fg})$ is compact by Corollary \ref{C:essSpectra} again. This implies the inclusion $\mathfrak{CT}\subseteq\mathcal{K}$, which completes the proof of statement (1).

For the proof of statement (2), consider the map $\Phi:$ $f\mapsto T_{f}+\mathcal{K}$ from $C(\bar{\mathbb{B}})$ into the quotient algebra $\mathfrak{T}/\mathcal{K} = \mathfrak{T}/\mathfrak{CT}$. It is clear that $\Phi$ is a $*$-homomorphism of $C^{*}$-algebras. By a result from the theory of $C^{*}$-algebras, the range of $\Phi$ is a closed $C^{*}$-subalgebra. On the other hand, it follows from the definition of $\mathfrak{T}$ that the range of the $\Phi$ is dense. Hence, $\Phi$ is a surjective $*$-homomorphism and we have $\mathfrak{T} = \{T_{f}+K: f\in C(\bar{\mathbb{B}}) \text{ and } K\in\mathcal{K}\}.$

From Corollary \ref{C:essSpectra}, the kernel of $\Phi$ is the ideal $\{f\in C(\bar{\mathbb{B}}): f|_{\mathbb{S}}\equiv 0\}$. Since the quotient of $C(\bar{\mathbb{B}})$ by this ideal is naturally isometrically $*$-isomorphic to $C(\mathbb{S})$, $\Phi$ induces a $*$-isomorphism $\tilde{\Phi}: C(\mathbb{S})\longrightarrow\mathfrak{T}/\mathcal{K}.$ This shows that the sequence
$$0\rightarrow \mathcal{K}\xrightarrow{\iota}\mathfrak{T}\xrightarrow{\tilde{\Phi}^{-1}\circ\pi} C(\mathbb{S})\rightarrow 0$$
is exact. Here $\iota$ is the inclusion map and $\pi:\mathfrak{T}\longrightarrow\mathfrak{T}/\mathcal{K}$ is the projection map.
\end{proof}

The following corollary shows that any Fredholm operator in $\mathfrak{T}$ has index zero. This generalizes Guo and Zheng's result to higher dimensions.

\begin{corollary}
Let $A$ be a Fredholm operator in $\mathfrak{T}$. Then $\ind(A)=0$.
\end{corollary}

\begin{proof}
By Theorem \ref{T:ToeplitzAlgebra}, there is a function $f\in C(\bar{\mathbb{B}})$ and a compact operator $K_1$ so that $A=T_f+K_1$. By the remark after Theorem \ref{T:reduction}, there is a compact operator $K_2$ so that $T_{f}=W^{*}M_{f_{*}}W+K_2$, where $f_{*}$ is the restriction of $f$ on $\mathbb{S}$. Since $A$ is a Fredholm operator, $M_{f_{*}}$ is a Fredholm operator on $L^2(\mathbb{S})$ and these two operators have the same index. But $M_{f_{*}}$ is Fredholm if and only if $f_{*}$ does not vanish on $\mathbb{S}$. In this case, $M_{f_{*}}$ is invertible (with inverse $M_{1/f_{*}}$) hence its index is $0$. The proof of the corollary is thus completed.
\end{proof}

\section{Toeplitz Operators with General Symbols}

In this section we consider Toeplitz operators with more general symbols. We present some necessary conditions for the compactness of $T_{f}$, where $f$ is assumed to be bounded. As a result, we show that if $f$ is a bounded harmonic function on $\mathbb{B}$ and $T_f$ is compact, then $f$ is the zero function.

We begin with a well known result that Toeplitz operators whose symbols have zero limit at the boundary of the unit ball are compact. The proof is based on the boundedness of kernel functions on compact subsets.

\begin{proposition}
\label{P:cptMultiplication}
If $f$ is a bounded (not necessarily continuous) function on $\mathbb{B}$ so that $\lim_{|x|\uparrow 1}f(x)=0$, then $M_f|_{b^2_{\nu}}$ is compact. As a consequence, the Toeplitz operator $T_f$ is compact on $b^2_{\nu}$.
\end{proposition}

\begin{proof}
For any $0<r<1$, let $\mathbb{B}_r=\{x\in\mathbb{R}^2: |x|\leq r\}$ and let $f_{r}=f\chi_{\mathbb{B}_r}$. It follows from the hypothesis that $\|f-f_r\|_{\infty}\to 0$ as $r\uparrow 1$. Therefore, $\|M_f-M_{f_r}\|\to 0$ as $r\uparrow 1$.

We now show that $M_{f_r}$ is a Hilbert-Schmidt operator on $b^2_{\nu}$ for $0<r<1$. Let $e_0, e_1, \ldots$ be an orthonormal basis for $b^2_{\nu}$. We have
\begin{align}\label{Eqn:cptMultiplication}
\sum_{j=0}^{\infty}\|M_{f_r}e_j\|^2 & = \sum_{j=0}^{\infty}\int_{\mathbb{B}}|f_r(x)|^2|e_j(x)|^2\mathrm{d}\nu(x)\notag\\
& = \int_{\mathbb{B}}|f_r(x)|^2\sum_{j=0}^{\infty}|e_j(x)|^2\mathrm{d}\nu(x)\\
& = \int_{\mathbb{B}_r}|f(x)|^2R(x,x)\mathrm{d}\nu(x).\notag
\end{align}
The last equality follows from the well known formula for the reproducing kernel function: $$R(x,x)=\|R_x\|^2 = \sum_{j=0}^{\infty}|\langle R_x,e_j\rangle|^2=\sum_{j=0}^{\infty}|e_j(x)|^2.$$
Since $R(x,x)$ is bounded for $x$ in $\mathbb{B}_r$ by \eqref{Eqn:ineqkernel}, the last integral in \eqref{Eqn:cptMultiplication} is finite. This shows that the operator $M_{f_r}|_{b^2_{\nu}}$ is a Hilbert-Schmidt operator on $b^2_{\nu}$. Therefore, $M_f|_{b^2_{\nu}}$, which is the norm limit of a net of compact operators, is compact. Since $T_f = PM_f|_{b^2_{\nu}}$, $T_f$ is also compact.
\end{proof}

The following proposition offers a necessary condition for a Toeplitz operator to be compact. 

\begin{proposition}
\label{P:cptToeplitz}
For each integer $m\geq 0$, let $\varphi_{m}$ be a positive function of the form $\varphi_{m}(x) = \sum_{j=1}^{s_m}|a_{j}^{(m)}(x)|^2,$ where $s_m$ is a positive integer and $a_j^{(m)}\in\mathcal{H}_m$ for $1\leq j\leq s_m$. Suppose $f$ is a bounded function on $\mathbb{B}$ so that $T_f$ is a compact operator on $b^2_{\nu}$. Then we have
\begin{equation}\label{Eqn:cptLimita}
\lim_{m\to\infty}\dfrac{\int_{\mathbb{B}}f(x)\varphi_{m}(x)\mathrm{d}\nu(x)}{\int_{\mathbb{B}}\varphi_m(x)\mathrm{d}\nu(x)}=0.
\end{equation}
In particular,
\begin{equation}\label{Eqn:cptLimitb}
\lim_{m\to\infty}\dfrac{\int_{[0,1)}\big(\int_{\mathbb{S}}f(r\zeta)\mathrm{d}\sigma(\zeta)\big)r^{2m}\mathrm{d}\mu(r)}{\int_{[0,1)}r^{2m}\mathrm{d}\mu(r)}=0.
\end{equation} 
\end{proposition}

\begin{proof}
Let $\epsilon>0$ be given. By Lemma \ref{L:cpt}, there is an integer $m_{\epsilon}>0$ so that for all integers $m\geq m_{\epsilon}$ and $p\in\mathcal{H}_{m}$, we have $\|T_fp\|\leq\epsilon\|p\|$. In particular, $\|T_fa_j^{(m)}\|\leq\epsilon\|a_j^{(m)}\|$ for $1\leq j\leq s_m$ for each such $m$. Now,
\begin{align*}
\Big|\int_{\mathbb{B}}f(x)\varphi_m(x)\mathrm{d}\nu(x)\Big| & = \Big|\sum_{j=1}^{s_m}\int_{\mathbb{B}}f(x)a_j^{(m)}(x)\bar{a}_{j}^{(m)}(x)\mathrm{d}\nu(x)\Big| \\
& = \Big|\sum_{j=1}^{s_m}\langle T_fa_j^{(m)},a_j^{(m)}\rangle\Big|\\
& \leq\sum_{j=1}^{s_m}\|T_fa_j^{(m)}\|\|a_j^{(m)}\|\\
& \leq\sum_{j=1}^{s_m}\epsilon\|a_j^{(m)}\|^2
 = \epsilon\int_{\mathbb{B}}\varphi_m(x)\mathrm{d}\nu(x).
\end{align*}
Therefore, $$\dfrac{|\int_{\mathbb{B}}f(x)\varphi_m(x)\mathrm{d}\nu(x)|}{\int_{\mathbb{B}}\varphi_m(x)\mathrm{d}\nu(x)}\leq\epsilon$$ for all $m\geq m_{\epsilon}$ and hence \eqref{Eqn:cptLimita} follows.

Let $e^{(m)}_1, \ldots, e^{(m)}_{h_m}$ be an orthonormal basis for $\mathcal{H}_m$ and define $$R_m(x,y)=\sum_{j=1}^{h_m}e^{(m)}_j(x)\bar{e}^{(m)}_{j}(y)$$ for $x,y\in\mathbb{B}$. Then $R_m$ is the reproducing kernel for $\mathcal{H}_m$ and since $\mathcal{H}_{m}$ is invariant under the action of the group of orthogonal transformations $O(n)$, $R_m(Ty,x)=R_m(y,T^{-1}x)$ for all $T\in O(n)$ and all $x, y$ in $\mathbb{B}$. The proof of these assertions is the same as the proof of \cite[Proposition 5.27]{AxlerSpringer2001}.

For $x\in\mathbb{B}$ and $T\in O(n)$, we have $R_m(Tx,Tx)=R_m(x,x)$. This implies that $R_m(x,x)=d_m|x|^{2m}$ for some constant $d_m$, which shows that $|x|^{2m} = d_m^{-1}\sum_{j=1}^{h_m}|e_j^{(m)}(x)|^2$. By choosing $\varphi_m(x)=|x|^{2m}$ and using integration in polar coordinates, we obtain
\begin{align*}
	\int_{\mathbb{B}}f(x)\varphi_m(x)\mathrm{d}\nu(x) & =\int_{\mathbb{B}}f(x)|x|^{2m}\mathrm{d}\nu(x) = \int_{[0,1)}\Big(\int_{\mathbb{S}}f(r\zeta)\mathrm{d}\sigma(\zeta)\Big)r^{2m}\mathrm{d}\mu(r),\\
	\int_{\mathbb{B}}\varphi_m(x)\mathrm{d}\nu(x)& =\int_{\mathbb{B}}|x|^{2m}\mathrm{d}\nu(x) = \int_{[0,1)}r^{2m}\mathrm{d}\mu(x).
\end{align*}
The limit \eqref{Eqn:cptLimitb} now follows from \eqref{Eqn:cptLimita}.
\end{proof}

\begin{theorem}\label{T:cptToeplitz}
Suppose $f$ is a bounded function on $\mathbb{B}$ so that the radial limit $f_{*}(\zeta)=\lim_{r\uparrow 1}f(r\zeta)$ exists for $\sigma$-almost every $\zeta$ on $\mathbb{S}$. If $T_f$ is a compact operator on $b^2_{\nu}$, then $f_{*}(\zeta)=0$ for $\sigma$-almost every $\zeta\in\mathbb{S}$.
\end{theorem}

\begin{proof}
For any multi-index $\alpha=(\alpha_1,\ldots,\alpha_n)\in\mathbb{N}_{0}^n$, the compactness of $T_f$ together with Corollary \ref{C:essSpectra} shows that $T_{f(x)x^{\alpha}}$ is a compact operator (here $x^{\alpha}=x_1^{\alpha_1}\cdots x_n^{\alpha_n}$). By Proposition \ref{P:cptToeplitz}, we have $$\lim\limits_{m\to\infty}\dfrac{\int_{[0,1)}\big(\int_{\mathbb{S}}f(r\zeta)r^{\alpha_1+\cdots+\alpha_n}\zeta^{\alpha}\mathrm{d}\sigma(\zeta)\big)r^{2m}\mathrm{d}\mu(r)}{\int_{[0,1)}r^{2m}\mathrm{d}\mu(r)}=0.$$
On the other hand, the hypothesis of the theorem together with the Dominated Convergence Theorem gives $$\lim_{r\uparrow 1} \int_{\mathbb{S}}f(r\zeta)r^{\alpha_1+\cdots+\alpha_n}\zeta^{\alpha}\mathrm{d}\sigma(\zeta)= \int_{\mathbb{S}}f_{*}(\zeta)\zeta^{\alpha}\mathrm{d}\sigma(\zeta).$$
It now follows from Lemma \ref{L:limitIntegrals} that $\int_{\mathbb{S}}f_{*}(\zeta)\zeta^{\alpha}\mathrm{d}\sigma(\zeta)=0$. Since this holds for any multi-index $\alpha$, we conclude that $f_{*}(\zeta)=0$ for $\sigma$-almost every $\zeta$ in $\mathbb{S}$.
\end{proof}

We say that a bounded function $f$ defined on $\mathbb{B}$ has a \textit{uniform radial limit} if there exists a function $f_{*}$ on $\mathbb{S}$ such that $$\lim_{r\uparrow 1}\Big(\sup\{|f(r\zeta)-f_{*}(\zeta)|: \zeta\in\mathbb{S}\}\Big)=0.$$ The function $f_{*}$ will be called the uniform radial limit of $f$. It is clear that any function $f$ in $C(\bar{\mathbb{B}})$ has a uniform radial limit, namely $f|_{\mathbb{S}}$. 

On the other hand, if $f_{*}$ is a bounded function on $\mathbb{S}$ and we define 
\begin{equation}\label{Eqn:extension}
\varphi(x)=\begin{cases}
|x|f_{*}(\frac{x}{|x|}) & \text{ if } 0<|x|\leq 1,\\
0 & \text{ if } x=0,
\end{cases}
\end{equation}
then it can be checked that $f_{*}$ is the uniform radial limit of $\varphi$. Moreover, if $f_{*}$ is continuous on $\mathbb{S}$, then $\varphi$ is continuous on $\bar{\mathbb{B}}$.

Using Proposition \ref{P:cptMultiplication} we extend Corollary \ref{C:essSpectra} to functions with continuous uniform radial limits.

\begin{corollary}\label{C:strongEssSpectra}
Let $f$ be a bounded function on $\mathbb{B}$ with uniform radial limit $f_{*}$ on $\mathbb{S}$. Assume that $f_{*}$ is continuous on $\mathbb{S}$. Let $g$ be a bounded function on $\mathbb{B}$. Then 
$T_{gf}-T_gT_f$ and $T_{gf}-T_fT_g$ are compact operators, $\sigma_{e}(T_f)=f_{*}(\mathbb{S})$ and $\|T_f\|_{e}=\sup\{|f_{*}(\zeta)|: \zeta\in\mathbb{S}\}$. In particular, $T_f$ is compact if and only if $f_{*}$ vanishes on $\mathbb{S}$. 
\end{corollary}
\begin{proof}
Let $\varphi$ be defined as in \eqref{Eqn:extension}. Then $\varphi$ belongs to $C(\bar{\mathbb{B}})$ and we have $\lim_{|x|\uparrow 1}|f(x)-\varphi(x)|=0$. Proposition \ref{P:cptMultiplication} shows that $T_{f}-T_{\varphi}$ and $T_{gf}-T_{g\varphi}$ are compact. The conclusions now follow from Corollary \ref{C:essSpectra}, which says $T_{g\varphi}-T_gT_{\varphi}$ and $T_{g\varphi}-T_{\varphi}T_g$ are compact operators, $\sigma_{e}(T_{\varphi})=\varphi(\mathbb{S})=f_{*}(\mathbb{S})$ and $\|T_{\varphi}\|_{e}=\sup\{|\varphi(\zeta)|: \zeta\in\mathbb{S}\}=\sup\{|f_{*}(\zeta)|: \zeta\in\mathbb{S}\}$.
\end{proof}

For functions whose uniform radial limits may not be continuous, we are unable to decide the validity of all conclusions in Corollary \ref{C:strongEssSpectra} but we do obtain a characterization for compactness.
\begin{corollary}\label{C:cptRadialLimit}
Suppose $f$ is a bounded function on $\mathbb{B}$ with uniform radial limit $f_{*}$ on $\mathbb{S}$. Then $T_f$ is compact if and only if $f_{*}(\zeta)=0$ for $\sigma$-almost every $\zeta$ in $\mathbb{S}$.
\end{corollary}
\begin{proof}
The ``if'' part follows from Proposition \ref{P:cptMultiplication} and the ``only if'' part follows from Theorem \ref{T:cptToeplitz}.
\end{proof}

Our last result in the paper is the fact that there is no non-zero compact Toeplitz operator with harmonic symbol. This was proved earlier by Choe, Koo and Na in \cite{ChoeRMJM2010} with a different approach for the case $\nu$ the normalized Lebesgue measure.

\begin{corollary}
Suppose $f$ is a bounded harmonic function on $\mathbb{B}$ so that $T_f$ is a compact operator on $b^2_{\nu}$, then $f(x)=0$ for all $x\in\mathbb{B}$.
\end{corollary}
\begin{proof}
It is well known (see \cite[Theorems 6.13 and 6.39]{AxlerSpringer2001}) that the radial limit $f_{*}(\zeta)=\lim_{r\uparrow 1}f(r\zeta)$ exists for $\sigma$-almost every $\zeta$ on $\mathbb{S}$ and $f$ is the Poisson integral of $f_{*}$. Since $T_f$ is assumed to be compact, Theorem \ref{T:cptToeplitz} shows that $f_{*}(\zeta)=0$ for $\sigma$-almost every $\zeta$ in $\mathbb{S}$. Thus $f(x)=0$ for all $x$ in $\mathbb{B}$.
\end{proof}

\providecommand{\bysame}{\leavevmode\hbox to3em{\hrulefill}\thinspace}
\providecommand{\MR}{\relax\ifhmode\unskip\space\fi MR }
\providecommand{\MRhref}[2]{%
  \href{http://www.ams.org/mathscinet-getitem?mr=#1}{#2}
}
\providecommand{\href}[2]{#2}

\end{document}